\RequirePackage{etoolbox}
\csdef{input@path}{{style/}{graphics/}}
\documentclass[MSNbibl,number,citesort,dvips]{arxbj}

\aid{0}
\volume{21}
\issue{1}
\pubyear{2015}
\firstpage{360}
\lastpage{373}
\doi{10.3150/13-BEJ570} 

\makeatletter
\newcommand{\eqref}[1]{(\ref{#1})}
\newtheorem{theorem}{Theorem}[section]
\newproclaim{assumption}[theorem]{Assumption}

\newproclaim{definition}[theorem]{Definition}
\newtheorem{lemma}[theorem]{Lemma}

\newproclaim{remark}[theorem]{Remark}
\newproclaim{example}[theorem]{Example}
\newremark{Trajectorial}{Trajectorial inequalities}

\newcommand{\E}{\mathbb{E}}
\renewcommand{\P}{\mathbb{P}}
\newcommand{\N}{\mathbb{N}}
\newcommand{\R}{\mathbb{R}}

\makeatother

\begin{document}
\begin{frontmatter}

\title{Pathwise versions of the Burkholder--Davis--Gundy inequality}
\runtitle{Pathwise versions of the Burkholder--Davis--Gundy inequality}

\begin{aug}
\author[a]{\inits{M.}\fnms{Mathias} \snm{Beiglb\"ock}\thanksref{e1}\ead[label=e1,mark]{mathias.beiglboeck@univie.ac.at}} \and
\author[a]{\inits{P.}\fnms{Pietro} \snm{Siorpaes}\corref{}\thanksref{e2}\ead[label=e2,mark]{pietro.siorpaes@maths.ox.ac.uk}}
\address[a]{Faculty of Mathematics, University of Vienna,
Nordbergstrasse 15, 1090 Wien, Austria.\\ \printead{e1,e2}}
\end{aug}

\received{\smonth{6} \syear{2013}}

%
\begin{abstract}
We present a new proof of the Burkholder--Davis--Gundy inequalities for
$1\leq p<\infty$. The novelty of our method is that these
martingale inequalities are obtained as consequences of elementary
\emph{deterministic} counterparts. The latter have
a natural interpretation in terms of robust hedging.
\end{abstract}

%
\begin{keyword}
\kwd{Burkholder--Davis--Gundy}
\kwd{martingale inequalities}
\kwd{pathwise hedging}
\end{keyword}

\end{frontmatter}

\section{Introduction}
In this paper, we derive estimates which compare the running maximum of
a martingale with its quadratic variation.
Given real numbers $x_n, h_n$, $ n\in\N$ we write
\[
x^*_n:= \max_{k\leq n} | x_k |,\quad\quad
[x]_n:= x_0^2+\sum
_{k=0}^{n-1} (x_{k+1}-x_k)^2,\quad\quad
(h\cdot x)_n:= \sum_{k=0}^{n-1}
h_k(x_{k+1}-x_k) .
\]
We will derive pathwise versions of the famous Burkholder--Davis--Gundy
inequalities.

\begin{theorem}\label{BDGThm}
For $1\leq p <\infty$, there exist constants $a_p,b_p<\infty$ such
that the following holds: for every $N\in\N$ and every martingale
$(X_k)_{k=0}^N$\def\theequation{BDG}
%
\begin{equation}
\label{BDGIneq}
\E[X]_N^{p/2}\leq a_p\E
\bigl[ \bigl(X_N^*\bigr)^p \bigr], \quad\quad\E \bigl[
\bigl(X_N^*\bigr)^p \bigr] \leq b_p
\E[X]_N^{p/2}.
\end{equation}
\end{theorem}

For $p\in(1, \infty)$ this was established by Burkholder \cite
{Bu66}. Under additional assumptions, Burkholder and Gundy \cite
{BuGu70} obtain a version for $p\in(0,1]$, while the case $p=1$ of
\eqref{BDGIneq} without restrictions is due to Davis \cite{Da70}.\def\theequation{\arabic{section}.\arabic{equation}}\setcounter{equation}{0}

For a modern account see, for instance, \cite{ChTe03}.

\begin{Trajectorial*}
The novelty of this note is that the above martingale inequalities are
established as consequences of \emph{deterministic} counterparts. We
postpone the general statements and first state the \emph{trajectorial
version} of Davis' inequality.
\end{Trajectorial*}

\begin{theorem}\label{PathDavisThm}
Let $x_0, \ldots, x_N$ be real numbers and set\footnote{Throughout
this paper we use the convention $0/0=0$.} $h_n:= \frac{x_n}{\sqrt {[x]_n + (x_n^*)^2}}, n\leq N$. Then
%
\begin{equation}
\label{PathDavis} \sqrt{[x]}_N \leq3 x^*_N - (h\cdot
x)_N,\quad\quad x^*_N \leq6 \sqrt{[x]}_N + 2 (h\cdot
x)_N.
\end{equation}
\end{theorem}

While the proof of Theorem~\ref{PathDavisThm} is not trivial, we
emphasize that the inequalities in \eqref{PathDavis} are completely
elementary in nature.
The significance of the result lies in the fact that it implies Davis'
inequalities:
indeed, if $(X_n)_{n=0}^N$ is a martingale, we may apply \eqref
{PathDavis} to each trajectory of $X$ and obtain a bounded and adapted
process $H$. The decisive observation is that, by the martingale property,
%
\begin{equation}
\label{Int=0}   \E\bigl[ (H\cdot X)_N\bigr]= 0,
\end{equation}
so Davis' inequalities (with $a_1=3, b_1=6$) follow from \eqref
{PathDavis} by taking expectations.

We recall that the BDG inequalities also apply if $X=(X_t)_{t}$ is a
cadlag local martingale, and that this follows from a straightforward
limiting procedure. Moreover, the inequalities are considerably simpler
to prove for \emph{continuous} local martingales (see, for example,
\cite{RoWi00}); in this case, they also hold for $p\in(0,1)$, as
proved by Burkholder and Gundy \cite{BuGu70}.

The problem of finding the optimal values of the constants $a_p, b_p$
is delicate, and has been open for 47 years and counting; we refer to
Osekowski \cite{Os10} for a discussion of the current state of research.

The rest of this paper is organized as follows.
In Section~\ref{history}, we discuss the history of the pathwise
approach to martingale inequalities.
In Section~\ref{Heuristics for the pathwise hedging approach}, we
explain the intuition behind the hedging strategy $h=(h_k)_k$ used in
the pathwise version of Davis' inequality.
In Section~\ref{Davis inequality for continuous local martingales}, we
give a short proof of one Davis' inequality for continuous martingales;
notably, this argument leads to a better constant compared to the
previous literature (to the best of our knowledge).
In Section~\ref{Davis inequality}, we establish Theorem~\ref{PathDavisThm}.
In Section~\ref{Pathwise Burkholder-Gundy Inequality}, we use Theorem~\ref{PathDavisThm} to derive trajectorial versions of the
BDG-inequalities in the $p>1$ case;
these also lead to their corresponding classical probabilistic
counterpart, thus concluding a fully analytic derivation of Theorem~\ref{BDGThm}.

\section{History of the trajectorial approach}
\label{history}\setcounter{equation}{0}

The inspiration of the pathwise approach to martingale inequalities
used in this paper comes from mathematical finance, more specifically,
the theory of model-independent pricing. The starting point of the
field is the paper \cite{Ho98a} of Hobson, which introduces the idea
to study option-prices by means of \emph{semi-static hedging}; we
explain the concepts using the inequality
%
\begin{equation}
\label{SemiStatic} \sqrt{[x]}_N \leq3 x^*_N - (h\cdot
x)_N
\end{equation}
appearing in Theorem~\ref{PathDavisThm}.
If the process $x=(x_n)_{n=0}^{N}$ describes the price evolution of
a financial asset, the functions $\Phi(x)=\sqrt{[x]}_N$ and $ \Psi
(x)=3 x^*_N$ have the natural financial interpretation of being exotic
options; specifically, here $\Phi$ is an option on realized variance,
while $\Psi$ is a look-back option.
The seller of the option $\Phi$ pays the buyer the amount $\Phi
(x_0,\ldots,x_N)$ after the option's expiration at time $N$, and
$(h\cdot x)_N$ corresponds to the gains or losses accumulated while
trading in $x$ according to the portfolio $h=(h_k)_k$.

The decisive observation of Hobson is that inequalities of the type
\eqref{SemiStatic} can be used to derive \emph{robust bounds} on the
relation of the prices of $\Phi$ and $\Psi$: independently of the
market model, one should never trade the option $\Phi$ at a price
higher than the price of $\Psi$, since the payoff $\Phi$ can be \emph
{super-hedged} using the option $\Psi$ plus self-financing trading.
Here the \emph{hedge} $3 x^*_N - (h\cdot x)_N$ is designated
\emph{semi-static}: it is made up of a static part -- the option $3
x^*_N$ which is purchased at time $0$ and kept during the entire time
range -- plus a dynamic part which corresponds to the trading in the
underlying asset according to the strategy $h$.

Since the publication of \cite{Ho98a} a considerable amount of
literature on the topic has evolved (e.g.,
\cite{Ro93,BrHoRo01a,HoPe02,CoHoPe08,DaObRa10,CoOb11a,CoOb11b,CoWa11,HoNe12,HoKl12});
we refer in particular to the survey by Hobson \cite{Ho11} for a very
readable introduction to this area. The most important tool in
model-independent finance is the Skorokhod-embedding approach; an
extensive overview is given by Ob{\l}{\'o}j in \cite{Ob04}. Starting
with the papers \cite{GaHeTo11,BeHePe12} the field has also been
linked to the theory of optimal transport, leading to a formal
development of the connection between martingale theory and robust
hedging (\cite{DoSo12,AcBePeSc13,DoSo13}). A benefit for the theory
of martingale inequalities is the following guiding principle:

\emph{Every martingale inequality which compares expectations of two
functionals has a deterministic counterpart.}

In fact, recently Bouchard and Nutz \cite{BoNu13} have coined this
into a rigorous theorem in the discrete time setup, see also \cite{BeNu14}.

This idea served as a motivation to derive the Doob-maximal
inequalities from deterministic, discrete-time inequalities in \cite{AcBePeScTe12}.\hskip.2pt\footnote{Notably, much of the approach of \cite
{AcBePeScTe12} was already developed earlier by Ob{\l}{\'o}j and Yor
\cite{ObYo06}.} In the present article, we aim to extend the approach
to the case of the Burkholder--Davis--Gundy inequalities.

\section{Heuristics for the pathwise hedging approach}
\label{Heuristics for the pathwise hedging approach}\setcounter{equation}{0}

The aim of this section is to explain the basic intuition which lies
behind the choice of the integrand in the pathwise Davis inequalities.
Arguments are simpler in the case of Brownian motion, which we will now
consider.

We focus on one of the two inequalities; according to the pathwise
hedging approach, we should be looking for a strategy $H$ and a
constant $a$ such that
$\sqrt{t} \leq a B^*_t + (H\cdot B)_t$.
Indeed, a reasonable ansatz to find a super hedging strategy is to
search for a function $f(b, b^*,t)$ such that
%
\begin{equation}
\label{BDG1Ansatz} \sqrt{t} \leq a B_t^* + \bigl(f\bigl(B, B^*,t\bigr)
\cdot B \bigr)_t,\quad\quad t\geq0 .
\end{equation}
To make an educated guess for the function $f$, we argue on a \emph
{purely heuristic} level and consider paths which evolve in a very
particular way. Assume first that the path $(B_t(\omega)_t)_{t\geq0}$
stays infinitesimally close to the value $b$ for all $t\geq t_0$:
we picture BM as a random walk on a time grid with size $\mathrm{d}t$, making
alternating up and down steps of height $\sqrt{\mathrm{d}t}$.
Thus, we assume that
$B$ evolves in the form
%
\begin{equation}
\label{DownUp} B_{t_0+2n \mathrm{d}t} = b,\quad\quad B_{t_0+(2n+1)
\,\mathrm{d}t}= b+\sqrt{\mathrm{d}t},\quad\quad n\geq0,
\end{equation}
where necessarily $b$ lies between $-B_{t_0}^*$ and $B_{t_0}^*$. The
left-hand side of \eqref{BDG1Ansatz} is of course increasing, so we
have to ensure the same behavior on the right side. A little
calculation reveals that this means that $f$ should have the form
%
\begin{equation}
\label{heu1} f\bigl( B, B^*, t\bigr)\approx-\frac{ B}{\sqrt t} \quad\quad\mbox{as } t\to
\infty;
\end{equation}
to see this, set $H_t:=f(B_t, B_t^*,t)$ and compare the value $\sqrt {t+2\,\mathrm{d}t}-\sqrt{t} \approx \mathrm{d}t/\sqrt{t} $
with
\begin{eqnarray*}
(H\cdot B)_{t+2\,\mathrm{d}t}-(H\cdot B)_t & \approx& f\bigl(t,b,b^*
\bigr) \,\mathrm{d}B_t+f \bigl(t+\mathrm{d}t,b+\sqrt{\mathrm{d}t},b^* \bigr)\, \mathrm{d}B_{t+\mathrm{d}t}
\\
& \approx& f\bigl(t,b,b^*\bigr) \sqrt{\mathrm{d}t}+f \bigl(t+\mathrm{d}t,b+\sqrt{\mathrm{d}t},b^* \bigr) (-
\sqrt{\mathrm{d}t})
\\
& \approx&- \bigl[f \bigl(t,b+\sqrt{\mathrm{d}t}, b^* \bigr)-f\bigl(t,b, b^*\bigr) \bigr]
\sqrt{\mathrm{d}t} + \mathrm{O}\bigl(\mathrm{d}t^{3/2}\bigr)
\\
& \approx&- f_b \,\mathrm{d}t.
\end{eqnarray*}
To assure that both sides of \eqref{BDG1Ansatz} grow at the same
speed, we thus need to require $\mathrm{d}t/\sqrt{t}\approx-f_b \,\mathrm{d}t$ which
leads to \eqref{heu1}.

Next, we consider a path which exhibits a different kind of extreme
evolution: assume that $B_t(\omega)\approx M t$ for some number $M>0$.
Simply setting $f( B, B^*, t)\approx- B/\sqrt t$ would lead to $
(f(B, B^*,t) \cdot B )_t \approx- 2M^2 t^{3/2}/3$. Taking $t$
sufficiently large, this quantity would eventually supersede $aB_t^*
\approx a M t$ independent of the choice of $a$, and thus \eqref
{BDG1Ansatz} would fail. So, this argument suggest to choose a function
which is bounded (at least for fixed $(t,B^*)$).
Moreover, dealing with a bounded integrand would conveniently allow to
follow the explanation after Theorem~\ref{PathDavis} and obtain Davis'
inequalities from the pathwise Davis' inequalities.
Thus, we could consider the function
%
\begin{equation}
\label{ContStrat} f\bigl(B,B^*, t\bigr)=-\frac{B_t}{\sqrt t\vee B^*_t}.
\end{equation}
Thanks to the additional term $a B_t^*$ in \eqref{BDG1Ansatz}, it is
not a problem if $f( B, B^*, t)\approx-2 B/\sqrt t$ is violated for
``small'' values of $t$;
and, if $\sqrt{t}$ is large compared to $B^*$, $f(B,B^*, t)\approx
-2B_t/\sqrt t$ holds, thus satisfying \eqref{heu1}.
Another similar possibility would be to use the function
%
\begin{equation}
\label{ContStrat2} f\bigl(B,B^*, t\bigr)=-\frac{B_t}{\sqrt{ t + (B^*_t)^2}},
\end{equation}
as in Theorem~\ref{PathDavisThm}; the latter turns out to lead to
easier computations in the discrete time case.
We choose however $f$ given by \eqref{ContStrat}
when dealing with continuous martingales, since this allows us to
obtain Davis' inequality with a better constant than the values we
could find in the literature.

\section{Davis inequality for continuous local martingales}
\label{Davis inequality for continuous local martingales}\setcounter{equation}{0}
We now derive one pathwise Davis' inequality for continuous local
martingales; integrating it yields the corresponding standard Davis' inequality.
We notice that Theorem~\ref{BDG1bt} provides the constant $3/2$, which
is smaller than the optimal constant for general cadlag martingales
(which is known to be $\sqrt{3}$, see \cite{bur02best}).
We do not address here the opposite pathwise Davis'\vadjust{\goodbreak} inequality
for continuous local martingales since we are only interested in
Theorem \ref{thm4.1} for illustrative purposes (since, as mentioned before,
Davis inequality for cadlag local
martingales follows from the case of martingales in discrete time).\vspace*{-1pt}

\begin{theorem}\label{thm4.1} If $M$ is a continuous local martingale such that
$M_0=0$ then\vspace*{-1pt}
%
\begin{equation}
\label{BDG1bt} \sqrt{[M]_t} \leq\frac{3}{2} M_t^*-
\biggl(\frac{M_t}{\sqrt{ [M]_t }
\vee M^*_t}\cdot M_t \biggr)_t \quad\quad\mbox{for all
} t\geq0 .
\end{equation}
\end{theorem}

\begin{pf}
By the Dambis--Dubins--Schwarz time change result, it is enough to
consider the case where $M$ is a Brownian Motion, which we will denote
by $B$.
From Ito's formula applied to the semi-martingales $B^2_t$ and $\sqrt
t\vee B_t^*$ we find\vspace*{-1pt}
\[
\mathrm{d}\frac{B_t^2}{\sqrt t\vee B_t^*}=-\frac{ B_t^2}{t \vee B_t^{*2} } \,\mathrm{d} \bigl(\sqrt t\vee B^*_t
\bigr) + \frac{1}{\sqrt t\vee B^*_t} ( 2B_t\, \mathrm{d}B_t+ \mathrm{d}t ).
\]
We may thus replace the integral in \eqref{BDG1bt} and arrive at the
equivalent formulation\vspace*{-1pt}
%
\begin{equation}
\label{BDG1tbE} \frac{B_t^2}{\sqrt t\vee B_t^*}+ \int_0^t
\frac{ B_s^2}{s \vee
B_s^{*2} }\, \mathrm{d} \bigl(\sqrt s\vee B^*_s \bigr) - \int
_0^t \frac{1}{\sqrt
s\vee B^*_s} \,\mathrm{d}s \leq3
B_t^*-2\sqrt t.
\end{equation}
Inequality \eqref{BDG1tbE} gets stronger if we replace each occurrence
of $B$ by $B^*$; thus, setting $f(t)= \sqrt t, g(t)= B_t^*$, it is
enough to prove the following claim:

Let $f, g\dvtx \R^+\to\R^+$ be continuous increasing functions such that
$f(0)=g(0)=0$ and $(f\vee g)(a)>0$ if $a>0$. Then, for all $a>0$\vspace*{-1pt}
%
\begin{equation}
\label{TechCore} \biggl(\frac{g^2}{f\vee g} \biggr) (a)+ \int_0^a
\frac{g^2}{f^2\vee
g^2} \,\mathrm{d}(f\vee g) - \int_0^a
\frac{1}{f\vee g} \,\mathrm{d}f^2\leq(3 g- 2 f) (a).
\end{equation}
To show this, observe that, by a change of variables
$\int\frac{g^2}{f^2\vee g^2} \,\mathrm{d}(f\vee g)= -\int g^2 \,\mathrm{d}\frac{f\vee
g}{f^2\vee g^2}$. Hence, integrating by parts on the interval
$(\varepsilon,a)$
and taking the limit $\varepsilon\to0$, we see that the left hand
side of \eqref{TechCore} equals\vspace*{-1pt}
\[
\int_0^a \frac{\mathrm{d}g^2-\mathrm{d}f^2}{f\vee g }.
\]
By a change of variables and applying trivial inequalities we obtain
\[
\int_0^a \frac{ \mathrm{d}g^2}{f\vee g}= \int
_0^a 1_{\{g>0\}} \frac
{\mathrm{d}g^2}{f\vee g} \leq
\int_0^a \frac{ 1_{\{g>0\}} \,\mathrm{d}g^2}{g} =2g(a) ,\quad\quad \int
_0^a \frac{\mathrm{d}f^2}{f\vee g} \geq \int
_0^a \frac{\mathrm{d}f^2}{f(\cdot)\vee g(a)} .
\]
If $f(a)\leq g(a)$, the last integral equals $f^2(a)/g(a)$; otherwise
there exists some $b\in[0,a)$ such that $f(b)=g(a)$, and then
evaluating separately the integral on $(0,b)$ and on $[b,a)$ we obtain that
\[
\int_0^a \frac{\mathrm{d}f^2}{f(\cdot) \vee g(a)}=
\frac{f^2(b)}{g(a)} +2\bigl(f(a)-f(b)\bigr)= 2f(a)-g(a) .
\]
Since $2y-x^2/y\leq3y-2x$ holds for $y>0$, either way \eqref
{TechCore} follows.\vadjust{\goodbreak}
\end{pf}

\section{Davis inequality}
\label{Davis inequality}\setcounter{equation}{0}

In this section, we prove Theorem~\ref{PathDavisThm}; in fact, we will
establish
that\footnote{Inequality (\ref{1half}) slightly improves on Inequality
(\ref{PathDavis}) by replacing the constant $3$ with the smaller
$1+\sqrt{2}$.}
%
\begin{eqnarray}
\label{1half} \sqrt{[x]_n} & \leq&(\sqrt{2}+1) x^*_n+(-h
\cdot x)_n,
\\
\label{2half} x^*_n & \leq&6 \sqrt{[x]_n}+(2h \cdot
x)_n ,
\end{eqnarray}
where the dynamic hedging strategy is defined by $ h_n=\frac
{x_n}{\sqrt{[x]_n + (x_n^*)^2}}$ as in Theorem~\ref{PathDavisThm}.

To prove \eqref{1half}, \eqref{2half} we introduce the convention,
used throughout the paper, that any sequence $(y_i)_{i\geq0}$ is
defined to be $0$ at time $i=-1$, and we define the auxiliary functions
$f,g$ for $m>0, q\geq0$, $|x|\leq m$ by\vspace*{-1pt}
%
\begin{eqnarray}
\label{firstf} f(x,m,q)&:=& -2\sqrt{q}+ \sqrt{m^2+q}-
\frac{m^2-x^2}{2\sqrt{m^2+q}},
\\
\label{firstg} g(x,m,q)&:=& -2m+ \sqrt{m^2+q}+\frac{m^2-x^2}{2\sqrt{m^2+q}}
\end{eqnarray}
and continuously extend them to $(x,m,q)=(0,0,0)$ by setting
$f(0,0,0)=g(0,0,0)=0$.
We will need the following lemma, whose proof is a somewhat tedious
exercise in calculus.\vspace*{-1pt}

\begin{lemma}
\label{DavisIneqfn}
For $d\in\R, |x| \leq m, q\geq0, m\geq0$ we have, with $c=\sqrt{2}-1$,
%
\begin{eqnarray}
\label{ImportantIneq2} f\bigl(x+d, m \vee|x+d| , q+d^2\bigr) - f(x,m,q) &
\leq&\frac{xd}{\sqrt {m^2+q}}+ \bigl( \sqrt{q+d^2} -\sqrt{q} \bigr),
\\[-0.5pt]
\label{ImportantIneq1} g\bigl(x+d, m \vee|x+d| , q+d^2\bigr) - g(x,m,q) &
\leq&- \frac{xd}{\sqrt {m^2+q}}+ c \bigl( \bigl(m\vee|x+d|\bigr) -m \bigr).\quad\quad
\end{eqnarray}
\end{lemma}

Before proving Lemma~\ref{DavisIneqfn} we explain why it implies
\eqref{1half} and \eqref{2half}.\vspace*{-1pt}

\begin{pf*}{Proof of Theorem~\ref{PathDavisThm}}
Since $f(x_0,|x_0|,x_0^2)\leq0$, \eqref{ImportantIneq2} gives\vspace*{-1pt}
\begin{eqnarray*}
-2 \sqrt{[x]_n} + x_n^*/2 &\leq& f\bigl(x_n,
x_n^*, [x]_n\bigr) \leq\sum_{k=0}^{n-1}
f\bigl(x_{k+1}, x^*_{k+1}, [x]_{k+1}\bigr)-f
\bigl(x_k, x_k^*, [x]_k\bigr)
\\[-0.5pt]
&\leq&(h\cdot
x)_n + \sqrt{[x]_n},
\end{eqnarray*}
which implies \eqref{1half}; and since $g(x_0,|x_0|,x_0^2)\leq0$, we
get \eqref{2half} from \eqref{ImportantIneq1} as follows\vspace*{-1pt}
\begin{eqnarray*}
-2 x_n^* + \sqrt{[x]_n} &\leq& g\bigl(x_n,
x_n^*, [x]_n\bigr) \leq\sum_{k=0}^{n-1}
g\bigl(x_{k+1}, x^*_{k+1}, [x]_{k+1}\bigr)-g
\bigl(x_k, x_k^*, [x]_k\bigr)
\\[-0.5pt]
&\leq&-(h\cdot
x)_n + c x^*_n.\vadjust{\goodbreak}
\end{eqnarray*}
\upqed\end{pf*}

Now we prove Lemma~\ref{DavisIneqfn}.

\begin{pf*}{Proof of Inequality \eqref{ImportantIneq1}}
It is enough to consider the case $m>0$, as the one where $m=0$ then
follows by continuity. Then, we can assume that $m=1$ through
normalization. Define $h(x,q,d)$ to be the LHS minus the RHS of \eqref
{ImportantIneq1}; since $h(x,q,d)=h(-x,q,-d)$, it is sufficient to deal
with the case $d\geq0$.

\noindent\textit{Case} I $[1\geq|x+d|]$:
Here we have to show that
%
\begin{equation}
\label{D1} h= \sqrt{1+q+d^2}+ \frac{1-(x+d)^2}{2\sqrt{1+q+d^2}}-\sqrt {1+q}-
\frac{1-x^2}{2\sqrt{1+q}}+\frac{xd}{\sqrt{1+q}}\leq0.
\end{equation}
Since $h_{xx} \geq0$, $h$ is convex, so it is sufficient to treat the
boundary cases $x=-1$ and $x=1-d$. To simplify notation, we set $r=
\sqrt{1+q}$; notice that $r\geq1$ and $0\leq d\leq2$.

\noindent\textit{Sub-case} I.A $[1\geq|x+d|, x=-1]$:
Then \eqref{D1} follows from
\begin{eqnarray*}
&&\sqrt{r^2+d^2}+\frac{1-(d-1)^2}{2\sqrt{r^2+d^2}}-r-\frac{d}r
\leq  0
\\
&&\quad\Leftarrow\quad r^2+d^2+d-d^2/2  \leq   (r+d/r)
\sqrt{r^2+d^2}
\\
&&\quad\Leftarrow\quad r^4+d^4/4 +d^2
+r^2d^2+d^3+2dr^2   \leq
r^4 +2dr^2+d^2+r^2d^2+2d^3+d^4
/r^2
\\
&&\quad\Leftarrow\quad d^4/4  \leq  d^3+ d^4/r^2,
\end{eqnarray*}
which is true since $0\leq d\leq2$.

\noindent\textit{Sub-case} I.B $[1\geq|x+d|, x=1-d]$: Here \eqref
{D1} amounts to
\begin{eqnarray*}
&&\sqrt{r^2+d^2}-r-\frac{1-(1-d)^2}{2r} +\frac{(1-d)d}{r}
 \leq  0
\\
&&\quad\Leftarrow\quad\sqrt{r^2+d^2}   \leq  r + d^2/2r
\\
&&\quad\Leftarrow\quad r^2+d^2   \leq  r^2 +
d^2 + d^4/4r^2.
\end{eqnarray*}

\noindent\textit{Case} II $[1\leq|x+d|]$:
Since $|x| \leq1$ and $d\geq0$, we find that $|x+d|\geq1$ implies
$x+d=|x+d| \geq1$. In this case $h$ equals
%
\begin{equation}
-(2+c) (x+d-1)+ \sqrt{(x+d)^2+q+d^2}-\sqrt{1+q}-
\frac{1-x^2}{2\sqrt {1+q}} + \frac{xd}{\sqrt{1+q}}.
\end{equation}
Since $s\mapsto\sqrt{s^2+1}$ is convex, $h \leq0$ holds iff it holds
for all $x$ on the boundary. Moreover if $-1 \leq1-d=x \leq1$, then
we already know that $h \leq0$ from the corresponding sub-case $1 \geq
|x+d|$; so we only need to show that $h \leq0$ for $x=1, q,d\geq0$
and for $x=-1, q \geq0, d \geq2$, respectively.

\noindent\textit{Sub-case} II.A [$1\leq|x+d|, x=1$]: We have to show
that, for all $q,d\geq0$,
\[
h(1,q,d)= -(2+c)d + \sqrt{(1+d)^2+q+d^2}- \sqrt{1+q}+
\frac{d}{\sqrt {1+q}}\leq0.
\]
Since $(1+d)^2+d^2=2(d+1/2)^2+1/2$ and $s\mapsto\sqrt{1+s^2}$ is
convex, it follows that $h(1,q,d)$ is convex in $d$; hence, the
inequality has to be checked only for $d=0$ and for $d \to\infty$.
The first case is trivial, and in the latter, after dividing both sides
by $d$, we arrive at
$-(2+c)+\sqrt{2} + 1/\sqrt{1+q}\leq0$, which holds by our choice of
$c$ and the fact that $q\geq0$.

\noindent\textit{Sub-case} II.B [$1\leq x+d, x=-1$]:
We have to show that, for all $q \geq0, d \geq2$,
\[
h(-1,q,d)=-(2+c) (d-2)+ \sqrt{(-1+d)^2+q+d^2}-
\sqrt{1+q}- \frac
{d}{\sqrt{1+q}} \leq0.
\]
As above, by convexity in $d$ it suffices to consider the cases $d=2$
and $d \to\infty$. The first one amounts to $\sqrt{5+q}\leq\sqrt {1+q} + 2/\sqrt{1+q} $, which is easily proved taking the squares. The
second one, after dividing by $d$, amounts to $-(2+c)+ \sqrt {2}-1/\sqrt{1+q}\leq0$,
which holds since $-(2+c)+ \sqrt{2}\leq0 $ by our choice of $c$.
\end{pf*}

\begin{pf*}{Proof of Inequality \eqref{ImportantIneq2}}
As before, we can assume w.l.o.g. that $m=1$ and $d\geq0$. Define
$k(x,q,d)$ to be the LHS minus the RHS of \eqref{ImportantIneq2}.

\noindent
\textit{Case} I [$1\geq|x+d|$]: In this case, $k$ equals
\[
\sqrt{1+q+d^2}- \frac{1-(x+d)^2}{2\sqrt{1+q+d^2}} -\sqrt{1+q}+\frac{1-x^2}{2\sqrt{1+q}} -
\frac{xd}{\sqrt{1+q}}- 3 \bigl( \sqrt{q+d^2} -\sqrt{q} \bigr) .
\]
Let us first isolate the terms that depend on $x$. Define
$k_0:=(1+q+d^2)^{-1/2}-(1+q)^{-1/2}$, and $k_2:=k- k_0(x+d)^2/2$, so that
\[
k_2=\sqrt{1+q+d^2} -\sqrt{1+q}- \frac{1}{2\sqrt{1+q+d^2}} +
\frac{1+d^2}{2\sqrt{1+q}} - 3 \bigl( \sqrt{q+d^2} -\sqrt{q} \bigr) .
\]
Notice that we can write
\[
k_0=\int_0^{d^2}
k_{1}(s)\, \mathrm{d}s \quad\quad\mbox{for } k_1(s):= \frac{\mathrm{d}}{\mathrm{d}s}
(1+q+s)^{-1/2},
\]
and similarly $k_2=\int_0^{d^2} k_{3}(s,d^2) \, \mathrm{d}s$ for
%
\begin{eqnarray}
  k_3\bigl(s,d^2\bigr)&:= &
\frac{\mathrm{d}}{\mathrm{d}s} \biggl( \sqrt{1+q+s} -\frac{1-s+d^2}{2\sqrt{1+q+s}}-3\sqrt{q+s} \biggr)
\\
&= & \frac{1}{2\sqrt{1+q+s}}+ \frac{2(1+q+s)+1-s+d^2}{4(1+q+s)^{3/2}} -\frac{3}{2\sqrt{q+s}}.
\end{eqnarray}
Since the $(k_i)_{i}$ do not depend on $x$ and $k_0\leq0$, $\max_{x}k
= k_2+k_0\min_{x} (x+d)^2/2$. Since $\min_{-1\leq x \leq1} (x+d)^2$
equals $0$ if $0\leq d \leq1$ and equals $(-1+d)^2$ if $1\leq d$, to
show $k\leq0$ we are lead to study the following two sub-cases.

\noindent\textit{Sub-case} I.A [$1\geq|x+d|, d\leq1$]: In this
case, $k=k_2$; to show that $k_2\leq0$ it is enough to show $k_3\leq0$.
Since $0\leq s\leq d^2\leq1$ we get $-s+d^2\leq1$, and so trivially
%
\begin{equation}
 k_3\leq\frac{2(1+q+s)+1+1}{4(1+q+s)^{3/2}}-\frac{2}{2\sqrt{q+s}}.
\end{equation}
So, calling $y:=q+s$, it is enough to prove that for all $y\geq0$
%
\begin{equation}
\label{poleq1} \frac{2y+4}{4(1+y)^{3/2}} - \frac{2}{2\sqrt{y}} \leq0,
\quad\quad\mbox{i.e.}\quad\quad
\sqrt{y}(y+2) \leq(1+y)^{3/2} 2,
\end{equation}
which is seen to be true by taking squares and bringing everything on
the RHS to obtain a polynomial whose coefficients are all positive.

\noindent\textit{Sub-case} I.B [$1\geq|x+d|, d\geq1$]:
In this case $k = k_2+k_0(1-d)^2/2$, so it is enough to show that
$k_3+k_1(1-d)^2/2 \leq0$. Since from $1\geq|x+d|, |x|\leq1$ it
follows that $d\leq2$, computations entirely similar\footnote{Use
that in this case $0\leq s\leq d^2\leq4$ implies $-s+d^2-(d-1)^2\leq
3$.} to the other sub-case
establish the desired result.

\noindent\textit{Case} II [$1\leq|x+d|$]: In this case,
$x+d=|x+d|\geq1$ and $k$ equals
\[
\sqrt{(x+d)^2+q+d^2} -\sqrt{1+q}+\frac{1-x^2}{2\sqrt{1+q}} -
\frac{xd}{\sqrt{1+q}}- 3 \bigl( \sqrt{q+d^2} -\sqrt{q} \bigr) .
\]
Since trivially $\mathrm{d}k/\mathrm{d}x\leq0$,
to show $k\leq0$ we can assume that $x=1-d$, in which case we can
write $k$ as $k=\int_0^{d^2} \tilde{k}(s) \,\mathrm{d}s$ for
%
\begin{eqnarray}
  \tilde{k}(s)&:= & \frac{\mathrm{d}}{\mathrm{d}s} \biggl( \sqrt{1+q+s} +
\frac{s}{2\sqrt {1+q}}-3\sqrt{q+s} \biggr)
\\
&= & \frac{1}{2\sqrt{1+q+s}}+ \frac{1}{2\sqrt{1+q}}-\frac{3}{2\sqrt{q+s}}.
\end{eqnarray}
Since $1-d=x\in[-1,1]$ we have $d^2\leq4$, and so to get $k\leq0$ it
suffices to show that $ \tilde{k} \leq0$ for $s\leq4$. This holds since
\[
\tilde{k}\leq\frac{1}{2\sqrt{1+q}}-\frac{2}{2\sqrt{q+s}} \leq0
\quad\quad\mbox{for } s\leq4 .
\]
\upqed\end{pf*}

\section{Pathwise Burkholder--Gundy inequality}
\label{Pathwise Burkholder-Gundy Inequality}\setcounter{equation}{0}

Garsia has given a simple proof of the fact that the
BDG inequalities for general $p\geq1$ are a consequence of
Davis inequality ($p=1$) and of the famous
lemma by Garsia and Neveu; in this section we revisit his proof and
turn it into pathwise discrete-time arguments.

Garsia's proof (for which we refer to \cite{Me76}, Chapter~3, Theorems 30 and
32 or to \cite{Ch75}) works similarly to how the Doob
$L^p$-inequalities for
$p>1$ follow by writing $x^p$ as an integral, applying the (weak) Doob
$L^1$-inequality, using Fubini's\vadjust{\goodbreak} theorem, and finally applying H\"
older's inequality (see for example \cite{ReYo99}).
The difference is that for the BDG
inequalities one needs to use a different integral expression for
$x^p$, and so one has to consider Davis' inequalities not on the time
interval $[0,T]$ but on $[\tau,T]$, where $\tau$ is a stopping time.

In the pathwise setting, by the guiding principle stated in Section~\ref{history}, if $L$ is a functional of a martingale $X$ and $\tau$
is a stopping time, a statement of the type $\E[L|F_\tau]\leq0$ will
have to be turned into one of the type $L+ (H\cdot X)_T- (H\cdot
X)_{\tau}\leq0$; moreover, since there will be no expectations
involved, H\"older's inequality will have to be replaced by Young's inequality.

We will need to consider discrete time stochastic integrals for which
the initial time is different from $0$; given $i<n$ and real numbers
$(h_j)_{i\leq j \leq n-1}$ and $(x_j)_{i\leq j \leq n}$, we define
%
\begin{equation}
\label{intfromi} (h\cdot x)_i^n:=\sum
_{j=i}^{n-1} h_j (x_{j+1}-x_j).
\end{equation}
Moreover if, for $i\leq j \leq n-1$, $h_j$ is a \emph{function} from
$\R^{j+1}$ to $\R$, given real numbers $(x_j)_{0\leq j \leq n}$ we define
$(h\cdot x)_i^n$ as
\[
\sum_{j=i}^{n-1} h_j(x_0,
\ldots, x_j) (x_{j+1}-x_j).
\]
Either way, we set $(h\cdot x)_i^n:=0$ if $n=i$.

We now deduce pathwise Davis' inequalities on $\{i,i+1,\ldots,n \}$
from the ones on $\{0,1,\ldots,n \}$ by a simple time shift.\vspace*{-1pt}

\begin{lemma}
\label{condDavis}
Assume that $\alpha,\beta>0$ and $h_n, k_n\dvtx  \R^{n+1} \to\R, n\geq
0$ satisfy
%
\begin{equation}
\label{1dav} \sqrt{[x]_n}   \leq\alpha x^*_n+(h \cdot
x)_n,\quad\quad x^*_n \leq\beta \sqrt{[x]_n} + (k\cdot
x)_n
\end{equation}
for every sequence $(x_n)_{n\geq0}$.
Define, for $i\geq0$, $n \geq i$, the functions $f^{(i)}_n, g^{(i)}_n\dvtx
\R^{n+1} \to\R$ by
\[
{f}^{(i)}_n\bigl((x_j)_{0\leq j\leq n}\bigr) :=
{h}_{n-i}\bigl((x_{l}-x_{i-1})_{i\leq l \leq n}\bigr),\quad\quad
{g}^{(i)}_n\bigl((x_j)_{j\leq
n}\bigr) :=
{k}_{n-i}\bigl((x_{l}-x_{i-1})_{i\leq l \leq n}\bigr).
\]
Then we have, for $n\geq i\geq0$,
\[
\sqrt{[x]_n}-\sqrt{[x]_{i-1}}  \leq2 \alpha
x^*_n + \bigl(f^{(i)}\cdot x\bigr)^n_i,\quad\quad
x^*_n -x^*_{i-1}\leq\beta\sqrt{[x]_n} +
\bigl(g^{(i)}\cdot x\bigr)^n_i.
\]
\end{lemma}

\begin{pf}Fix $n\geq i\geq0$, $(x_n)_{n\geq0}$ and let
$y^{(i)}_j:=x_{j+i}-x_{i-1}$. Applying \eqref{1dav} to
$(y^{(i)}_j)_{j\geq0}$ we find
\begin{eqnarray*}
\sqrt{[x]_n}-\sqrt{[x]_{i-1}}&\leq&\sqrt{[x]_n-[x]_{i-1}}
=\sqrt {\bigl[y^{(i)}\bigr]_{n-i}}   \leq \alpha
\bigl(y^{(i)}\bigr)^*_{n-i} +\bigl(h \cdot y^{(i)}
\bigr)_{n-i}
\\
&\leq&\alpha2 x^*_n+\bigl(f^{(i)} \cdot x
\bigr)_{i}^n,
\end{eqnarray*}
and (respectively)
\[
x^*_n - x_{i-1}^* \leq\bigl(y^{(i)}
\bigr)^*_{n-i}  \leq\beta\sqrt {\bigl[y^{(i)}
\bigr]_{n-i}}+\bigl(k \cdot y^{(i)}\bigr)_{n-i} \leq
\beta\sqrt {[x]_n}+\bigl(g^{(i)} \cdot x\bigr)_{i}^n
.
\]
\upqed\end{pf}

Here follows the pathwise version of Garsia--Neveu's lemma.

\begin{lemma}
\label{GarsiaNeveu}
Let $p>1$, $ c_n\in\R, (x_j)_{j\leq n}, (h^{(i)}_n)_{i\leq n} \in\R
^{n+1}$, and assume that $0=a_{-1} \leq a_0 \leq\cdots\leq a_n
<\infty$ and
\[
a_n -a_{i-1} \leq c_n + \bigl(h^{(i)}
\cdot x\bigr)_i^n \quad\quad\mbox{for } n\geq i\geq 0 .
\]
Then, if we set
\[
w_j:=\sum_{i=0}^j p
\bigl(a^{p-1}_i - a^{p-1}_{i-1} \bigr)
h^{(i)}_j,\quad\quad j\leq n ,
\]
we have that
%
\begin{eqnarray}
\label{gnlc1} a_n^p& \leq& pc_n
a_{n}^{p-1}+ (w\cdot x)_n ,
\\
\label{gnlc2} a_n^p& \leq&(p-1)^{p-1}
c_n^p + (pw\cdot x)_n .
\end{eqnarray}
\end{lemma}

\begin{pf}
From $a^p_n = p(p-1)\int_0^{a_n} s^{p-2} (a_n-s) \,\mathrm{d}s = p\sum_{i=0}^{n}
\int^{a_i}_{a_{i-1}} (p-1) s^{p-2} (a_n-s) \,\mathrm{d}s$ and $a_n-s\leq
a_n-a_{i-1}$ on $s\in[a_{i-1},a_i]$, we find
\eqref{gnlc1} by writing
\begin{eqnarray*}
a_n^p &\leq& p  \sum_{i=0}^{n}
\bigl(a^{p-1}_i - a^{p-1}_{i-1}\bigr)
(a_n-a_{i-1})
\\
&\leq& p  \sum_{i=0}^{n}
\bigl(a^{p-1}_i - a^{p-1}_{i-1}\bigr)
\bigl[c_n+ \bigl(h^{(i)} \cdot x \bigr)_i^n
\bigr]
\\
&= & pc_n a_{n}^{p-1} + p  \sum
_{i=0}^{n} \sum_{j=i}^{n-1}
\bigl(a^{p-1}_i - a^{p-1}_{i-1}\bigr)
h_j^{(i)} (x_{j+1}-x_j)
\\
&= & pc_n a_{n}^{p-1} +  \sum
_{j=0}^{n-1} \Biggl( \sum_{i=0}^{j}
p \bigl(a^{p-1}_i - a^{p-1}_{i-1}\bigr)
h_j^{(i)} \Biggr) (x_{j+1}-x_j)  =
pc_n a_{n}^{p-1} + (w\cdot x)_n.
\end{eqnarray*}
We then obtain (\ref{gnlc2}) from (\ref{gnlc1}) by
applying Young's inequality $ab\leq C_\epsilon a^p/p+\epsilon b^q/q$
(where $C_\epsilon^{-1}=p(\epsilon q)^{p-1}$ and $1/p+1/q=1$) with
$\epsilon=1/p$, $a=c_n, b=a_{n}^{p-1}$.
\end{pf}

Finally, from Theorem~\ref{PathDavisThm}, Lemma~\ref{condDavis} and
Lemma~\ref{GarsiaNeveu}, we obtain the following discrete-time
pathwise BDG inequalities for $p>1$.
We recall that, by convention, $x_{-1}=x^*_{-1}=[x]_{-1}=0$ and
$0/0=0$, and in particular the integrand $f_n^{(i)}$ is well defined.

\begin{theorem}\label{PathBGThm}
Let $x_0, \ldots, x_N$ be real numbers, $c_p:= 6^p(p-1)^{p-1}$ for
$p>1$, and define
\[
h_n:= \sum_{i=0}^n
p^2 \Bigl(\sqrt{[x]_i^{p-1}} -\sqrt
{[x]_{i-1}^{p-1}} \Bigr) f^{(i)}_n,\quad\quad
g_n:= \sum_{i=0}^n
p^2 \bigl(\bigl(x^*_i\bigr)^{p-1} -
\bigl(x^*_{i-1}\bigr)^{p-1} \bigr) f^{(i)}_n,
\]
where
\[
f^{(i)}_n:=\frac{x_n-x_{i-1}}{\sqrt{[x]_n-[x]_{i-1} + \max_{i\leq
k\leq n}(x_{k}-x_{i-1})^2}}.
\]
Then
%
\begin{equation}
\label{PathBG} \sqrt{[x]_N^{p}} \leq
c_p \bigl(x^*_N\bigr)^p - (h\cdot
x)_N,\quad\quad \bigl(x^*_N\bigr)^p \leq c_p
\sqrt{[x]_N^{p}} + 2(g\cdot x)_N.
\end{equation}
\end{theorem}

We notice that Theorem~\ref{PathBGThm}
yields \eqref{BDGIneq}; indeed, given a finite constant $N$ and a
martingale $(X_n)_{n=0}^N$, trivially $\sqrt{[X]_N}$ and $X_N^*$ are
in $L^p(\P)$ iff $X_n$ is in $L^p(\P)$ for every $n\leq N$, and in
this case the adapted integrands $(H_n)_{n=0}^{N-1}$ and
$(G_n)_{n=0}^{N-1}$ which we obtain applying Theorem~\ref{PathBGThm}
to the paths of $X$ are in $L^q(\P)$ for every $n$ (for $q=p/(p-1)$), thus
$H\cdot X$ and $G\cdot X$ are martingales and so
\[
  \E\bigl[ (H\cdot X)_N\bigr]= 0= \E\bigl[ (G\cdot
X)_N\bigr],
\]
and the Burkholder--Davis--Gundy inequalities for $p>1$ (with
$a_p=b_p=6^p(p-1)^{p-1}$) follow from \eqref{PathBG} by taking
expectations, completing the proof of Theorem~\ref{BDGThm}.

\section*{Acknowledgements}

The authors thank Harald Oberhauser for comments on an
earlier version of this paper. The first author thanks the Austrian
Science Fund for support through project p21209.

%

\printhistory

\end{document}